\documentclass[notitlepage,leqno,10pt]{article}

\newenvironment{proof}{\noindent{\sc Proof}.\enspace}{\rule{2mm}{2mm}\medskip}

\usepackage{amsmath}
\usepackage{amssymb}

\newtheorem{theorem}{Theorem}[section]
\newtheorem{lemma}[theorem]{Lemma}

\newcommand{\supp}{{ \rm supp}}

\title{Critical regularity for   elliptic equations
from Littlewood-Paley theory }

\author{Denis A. Labutin}

\date{}

\begin{document}

\maketitle

\begin{center}

{\small Department of mathematics, 
University of California, Santa Barbara,
CA 93106, USA}

\end{center}

\footnotetext[1]{E-mail address: 
labutin@math.ucsb.edu}

\

\

\noindent {\sc abstract}. Using simple facts from harmonic analysis, namely
Bernstein inequality  and  Plansherel isometry, 
we prove that the  pseudodifferential equation
$\Delta^\alpha u+Vu=0$
improves the  Sobolev regularity of solutions
provided the potential
$V$
is integrable with the critical power
$n/2\alpha>1$.

\begin{center}

\bigskip\bigskip

\centerline{\bf AMS subject classification:
35J60, 58J05, 58J40 }

\end{center}

%%%%%%%%%%%%%%%%%%%%%%%%%%%%%%%%%%%%%%%%%%%%%
%
%                    PAPER
%
%%%%%%%%%%%%%%%%%%%%%%%%%%%%%%%%%%%%%%%%%%%%%

%%%%%%%%%%%%%%%%%%%%%%%%%%%%%%%%%%%%%%%%%%%%
%
%  The command \large  is inserted 
%
%%%%%%%%%%%%%%%%%%%%%%%%%%%%%%%%%%%%%%%%%%%%

%\large

%
%
%
%
%
%
%
%      MAIN TEXT STARTS HERE
%
%
%
%
%
%
%
%
%

\section{Introduction}

\setcounter{equation}{0}
In this paper we prove the following local regularity  
result for (complex pseudodifferential) elliptic equations.
Let 
$B_R$
denote a ball  of radius 
$R$
and
center
$0$
in 
$\mathbf{R}^n$,
$n=2$,
$3$,
$\ldots$.
%
%
%
%%%%%%%%%%%%%%%%%%%%%%%%%%%%%%%%%%%%%%%%%%%%%%%%%%%
%
%                    MAIN THEOREM
%
%%%%%%%%%%%%%%%%%%%%%%%%%%%%%%%%%%%%%%%%%%%%%%%%%%%
%
%
\begin{theorem}
\label{maintheorem}
Let
$u\in H^s(B_1)$
solve
\begin{equation}
\label{maineq}
\Delta^\alpha
u
+
Vu
=0
\quad
{\it in}
\quad
B_1
\end{equation}
with
$V\in L^{n/2\alpha}(B_1)$,
$0<2 \alpha< n$,
and
$0<2s<n$.
Assume also that
\begin{equation}
\label{restriction_on_s}
2\alpha -\frac{n}{2} <s< 2\alpha.
\end{equation}
Then
\begin{equation}
\label{smoothness_improved}
u\in
H^{s+\varepsilon} (B_{1/2})
\end{equation}
for some 
$\varepsilon>0$,
$\varepsilon=\varepsilon(n, \alpha, s)$.
\end{theorem}

By 
$H^s$
we denote the Sobolev Hilbert space
of order
$s\in\mathbf{R}^1$.
We write
$u\in H^s(B_R)$
if
$u$ 
is 
a distribution on
$B_R$
such that 
\begin{equation*}
u=\tilde{u}|_{B_R}
\quad
{\rm
for \quad
some
}
\quad
\tilde{u}\in H^s(\mathbf{R}^n).
\end{equation*}
For any
$\alpha>0$
the
operator
$\Delta^\alpha$
is a nonlocal pseudodifferential operator
defined on 
$H^s(\mathbf{R}^n)$,
\begin{equation*}
\Delta^\alpha
\colon
H^s(\mathbf{R}^n)
\longrightarrow
H^{s-2\alpha}(\mathbf{R}^n),
\end{equation*}
for any
$s\in\mathbf{R}^{n}$.

Thus the theorem
states  that, if
a  distribution
$u\in H^s(B_1)$
has an extension
$\tilde{u}$
with
$\Delta^{\alpha} \tilde{u} \in H^{s-2\alpha}(\mathbf{R}^{n})$
satisfying
(\ref{maineq})
in
$B_1$,
then
we can find  a distribution from 
$H^{s+\varepsilon}(\mathbf{R}^n)$,
$\varepsilon>0$,
coinciding with 
$u$
in
$B_{1/2}$.
The particular case of integer
$s=\alpha$ 
arises in the calculus of variations. 
In this case our  theorem says that 
equation
(\ref{maineq})
improves 
the regularity
of 
$H^\alpha$-solutions 
if
$0<2\alpha<n$.

We will prove the theorem using only two
simple facts from Littlewood-Paley theory, namely,
the Plansherel isometry and the Bernstein inequality.
Regularity for 
(\ref{maineq})
for  the end-point 
relations between parameters
$n$,
$\alpha$,
and
$s$
cannot apparently be established with such 
simple tools. 
The only non-obvious assumption on the parameters in
Theorem~\ref{maintheorem}
is the lower bound for 
$s$
in
(\ref{restriction_on_s}).
Together with the Sobolev embedding it garantees that
the product
$Vu$
is defined as a distribution.

Theorem~\ref{maintheorem}
can be derived from the results of Y.Y.Li
\cite{YYLi},
see Theorem 1.3 there.
Equation
(\ref{maineq})
is treated  in
\cite{YYLi}
as an integral equation in the physical space 
and the frequency space is not used there at all.
Technique developed in the present paper
does not depend on the structure of the fundamenatal
solution of
$\Delta ^\alpha$. 
In particular,
it allows to establish the local regularity for more general
pseudodifferential equations on smooth manifolds
\cite{Labutin2}.

The function 
$V$
is integrable with the {\it critical} 
power in the theorem
meaning the following: 
if
$V\in L^p(B_1)$
with
$p<n/2\alpha$ 
then in general
(\ref{smoothness_improved})
does not hold for  any
$\varepsilon>0$
as the family of examples below shows.
If
$V\in L^p(B_1)$
with
$p>n/2\alpha$ 
then
the improved regularity
(\ref{smoothness_improved})
is easy to prove.
Indeed, in this case 
Sobolev and Holder inequalities 
imply at once that 
\begin{equation}
\label{supercritical_equation}
\Delta^\alpha
u
=f
\quad
{\rm in}
\quad
B_1
\quad
{\rm with}
\quad
f\in L^p(B_1),
\quad
\frac{1}{p}
<
\frac{2\alpha}{n}
+
\frac{1}{2}
-
\frac{s}{n}.
\end{equation}
Now
(\ref{smoothness_improved}) is
the straightforward consequence 
of the Calderon-Zygmund estimate
and Sobolev inequality.

The main purpose and application  of 
Theorem~\ref{maintheorem}
is deriving the 
{\it full regularity}
for quasilinear (complex pseudodifferential)
elliptic equations with the critical growth nonlinearity.
For such  application any
$\varepsilon>0$
in
(\ref{smoothness_improved})
works equally well. It is for this reason that we do not 
care about the sharp  value of
$\varepsilon(n, \alpha, s)$
in 
Theorem~\ref{maintheorem}.
For example, consider a weak solution
$u\in H^\alpha$,
$0<2\alpha<n$,
of the equation
\begin{equation*}
\Delta^\alpha u
+g(x,u)
=0
\quad
{\rm in}
\quad
B_1.
\end{equation*}
Assume that
$g$
is a smooth, possibly complex valued function of the critical growth:
\begin{equation*}
|g(x,t)|
\leq
C\left(
1+
|t|^{(n+2\alpha)/(n-2\alpha)}
\right)
\quad
{\rm for 
\quad 
all}
\quad
x\in B_1, 
\,
t\in\mathbf{C}^1.
\end{equation*}
We can write
\begin{eqnarray*}
g(x, u)
&
=
&
\frac{g(x, u)}{1+|u|}
+
\left(
\frac{g(x, u)}{1+|u|}
\,
\frac{|u|}{u}
\right)u
\\
&
=
&
f+Vu
\end{eqnarray*}
with
$f$
as in
(\ref{supercritical_equation})
and
$V\in L^{n/2\alpha}$.
Now the application of
Theorem~\ref{maintheorem}
combined with 
Calderon-Zygmund and Sobolev inequalities
improve the integrability of 
$u$.
Then 
Schauder estimates
imply that
$u\in C^\infty(B_{1/2})$.

This way of proving the
regularity for the critical semilinear equations 
was suggested  for
$\alpha =1$
by 
Brezis and Kato
\cite{Brezis_Kato},
see also
Appendix B
in
\cite{Struwe_book}.
These authors improved 
{\it integrability} 
of
$u$
using Moser's iteration technique.
Earlier Trudinger 
\cite{Trudinger_Pisa}
(also in  the case $\alpha =1$)
had already 
used Moser's iterations to 
prove the full regularity for
the nonlinear problem directly.
The case of an integer 
$\alpha >1$
has attracted recent attention  in
\cite{Chang_Gursky_Yang},
\cite{Uhlenbeck_Viaclovsky}
due to 
its applications in conformal geometry.
In a related  paper
\cite{YYLi}
Y. Y. Li
proved the full regularity for the equation
\begin{equation}
\label{y_y_li}
\Delta^\alpha
u
+
u^
{
(n+2\alpha)/(n-2\alpha)
}
=0,
\quad
u>0,
\quad
0<2\alpha<n.
\end{equation}
The main goal in
\cite{YYLi}
was to establish  Liouville-type theorems for
(\ref{y_y_li})
in
$\mathbf{R}^n$
using the moving spheres method.
Earlier
Liouville
theorems were
proved for
(\ref{y_y_li})
in 
\cite{Lin},
\cite{Chang_Yang},
\cite{Chen_Li_Ou}
using the moving plane method. 
The Littlewood-Paley approach
was used in
\cite{Chemin_Xu}
and in
\cite{Tao}
to give proofs  of regularity of Holder-continuous
harmonic maps and harmonic maps  from surfaces into spheres respectively.

Elliptic equations
with 
{\it supercritical} 
nonlinearity
do not improve the regularity of  solutions.
For example,
for 
$\alpha =1$
and
any
$p>(n+2)/(n-2)$,
$n\geq 3$,
the function
\begin{equation*}
u(x)=\frac{A}{|x|^a},
\quad
a=2/(p-1),
\quad
A=(a(n+a-2))^{1/(p-1)},
\end{equation*}
satisfies
\begin{equation*}
\Delta u
+u^p
=0
\quad
{\rm in}
\quad
B_1,
\quad
u\in H^1(B_1).
\end{equation*}
However,
$u$
is not smooth in
$B_{1/2}$.
Pohozaev in
\cite{Pohozaev}
investigated local regularity for 
supercritical semilinear  problems, and
established some sharp low regularity 
results.

\noindent
{\bf Acknowledgments.} 
This work was done when the author was visiting the  
Australian National University in 2003 by the invitation of 
Neil Trudinger and Xu-Jia Wang.  
The author also wishes to thank O. V. Besov, M. L. Goldman,  S. I. Pohozaev, and 
other participants of the 
Fall, 2004 seminar on function spaces at the Steklov Institute
for their comments.

\section{Proof of Theorem~\ref{maintheorem}}
\setcounter{equation}{0}

Let
$\{\widehat{\varphi}_j\}_{j=-\infty}^{+\infty}$
be the standard 
smooth partition of unity in the 
Littlewood-Paley theory
\cite{Petree},
\cite{Triebel},
\cite{Stein1},
\cite{Stein2}.
Thus 
$\widehat{\varphi}_j=\widehat{\varphi}(\cdot/ 2^j)$
is supported in, say,
the ring
\begin{equation*}
\{ \xi\colon 2^j 3/5 \leq |\xi| \leq 2^j 5/3\}
\subset
\left(
B_{2^{j+1}}\setminus B_{2^{j-1}}
\right).
\end{equation*}
Let 
$P_j$
denote
the Littlewood-Paley projection,
\begin{equation*}
(P_j f)^\wedge
=\widehat{\varphi}_j \widehat{f},
\quad
f\in
\mathcal{S}'.
\end{equation*}
We also set
\begin{equation*}
P_{a<\cdot<b}
=
\sum_{j=a+1}^{b-1}
P_j.
\end{equation*}
Distributions  with the localised Fourier transform 
enjoy the important Bernstein inequality:
for
$f\in\mathcal{S}'$
and
$1\leq p\leq q\leq\infty$
\begin{equation*}
\|f\|_q
\lesssim
2^{nj((1/p)-(1/q))}
\|f\|_p
\quad
{\rm provided}
\quad
\supp\,\widehat{f}
\subset
B_{2^j}.
\end{equation*}
For
$s\in\mathbf{R}^1$
the Sobolev space
$H^s(\mathbf{R}^n)$
consists of distributions with the finite norm
\begin{equation*}
\|f\|_{H^s}
=
\|P_{\cdot\leq 0} f\|_2
+
\left(
\sum_{j=1}^\infty
2^{2js}
\|P_j f\|_2^2
\right)^{1/2}.
\end{equation*}
The Plansherel isometry implies that for
$s=1$,
$2$,
$\ldots$
the space
$H^s(\mathbf{R}^n)$
consists of distributions  with all derivatives
up to the order
$s$
lying
in
$L^2(\mathbf{R}^n)$.

%\noindent
%{\it Proof.}
\begin{proof}{\sc
(of Theorem~\ref{maintheorem})}
{\bf 1.}
First,  we localise the problem.
Take a cutoff function
$\eta_\rho$,
\begin{equation*}
\eta_\rho =1
\quad
{\rm in}
\quad
B_\rho,
\quad
\eta_\rho =0
\quad
{\rm outside}
\quad
B_{2\rho}.
\end{equation*}
The commutator of
the multiplication by 
$\eta_{\rho}$
and 
$\Delta^\alpha$
is a pseudodifferential operator of order
$2\alpha -1$.
For interger
$\alpha$
this is just the Leibnitz formula for the derivative
of the product.
Hence, for some
$F\in H^{s-2\alpha +1}(\mathbf{R}^n)$
we obtain
\begin{eqnarray}
\Delta^\alpha(\eta_\rho u)
&
=
&
\Delta^\alpha
(\eta_\rho \tilde{u})
\nonumber
\\
&
=
&
\eta_\rho \Delta^\alpha \tilde{u}
+F
\nonumber
\\
&
=
&
-(\eta_{2\rho} V)
(\eta_\rho u)
+F
\quad
{\rm in}
\quad
\mathbf{R}^n.
\label{localised_equation}
\end{eqnarray}
To economize on notations denote
$u\eta_\rho$
by 
$u$
and
$V\eta_{2\rho}$
by
$V$.
Then in 
(\ref{localised_equation})
we have
$u\in H^s (\mathbf{R}^n)$,
$\supp(u) \subset B_{2\rho}$.
Moreover, 
the 
$L^{n/2\alpha}$-norm 
of
$V$
is small
when
$\rho$
is small.
In the proof we,  by
making
this norm
small enough, 
will establlish
that
\begin{equation*}
u
\,
(=\eta_\rho u)
\, 
\in 
H^{s +\varepsilon} (\mathbf{R}^n).
\end{equation*}
Statement
(\ref{smoothness_improved})
then follows by covering
$B_{1/2}$
with  small balls.
Therefore the goal is to 
choose a suitable
$\rho$
so 
that
for some constant
$C>0$,
$C=C(u,V, \rho, n,\alpha,s)$,
\begin{equation}
\label{goal_estimate}
\|P_k u\|_2
\leq
\frac{C}{2^{(s+\varepsilon)k}}
\quad
{\rm 
for
\quad 
all}
\quad
k\geq 1.
\end{equation}
Clearly it is enough to prove
(\ref{goal_estimate})
only for large
$k$.

{\bf 2. }
The product  in the right hand side of
(\ref{localised_equation})
is an integrable function as a result of
(\ref{restriction_on_s}).
Hence, applying the Littlewood-Paley projection,
we derive that
\begin{eqnarray}
\label{first_estimate}
2^{2\alpha k}
\|P_k u\|_2
&
\lesssim
&
\|P_k (Vu)\|_2
+
\|P_k F\|_2
\nonumber
\\
&
\lesssim
&
\|P_k(Vu)\|_2
+
C_F 2^{(2\alpha-s-1)k}.
\end{eqnarray}
Thus to prove
(\ref{goal_estimate})
we need to estimate
$P_k(Vu)$.
We take into account the 
localisation of the 
Littlewood-Paley projections 
in the frequency space. 
It implies  that
for 
$f,g\in\mathcal{S}'$
the distribution 
$P_k(P_i f P_j g)$
vanishes identically if
\begin{equation*}
\Big(
B_{2^{i+1}}\setminus B_{2^{i-1}}
+
B_{2^{j+1}}\setminus B_{2^{j-1}}
\Big)
\cap
\Big(
B_{2^{k+1}}\setminus B_{2^{k-1}}
\Big)
=\emptyset.
\end{equation*}
Consequently
for a fixed
$k\in\mathbf{Z}$
\begin{eqnarray}
\label{decomposition}
P_k(Vu)
&
=
&
\sum_{i,j\in \mathbf{Z}}
P_k(
P_i V
P_j u
)
\nonumber
\\
&
=
&
\left\{
\sum_{i,j\in LL}
+
\sum_{i,j\in LH}
+
\sum_{i,j\in HL}
+
\sum_{i,j\in HH}
\right\}
P_k(
P_i V
P_j u
)
\nonumber
\\
&=&
I+II+III+IV,
\end{eqnarray}
where
$LL$,
$LH$,
$HL$,
and 
$HH$
are the
low-low,
low-high,
high-low,
and
high-high
frequencies interaction zones
on the integer lattice:
\begin{eqnarray*}
LL
&=&
\left\{
i,j\in \mathbf{Z}
\colon
\
k-5\leq i,j\leq k+7,
\
\min\{i,j\}\leq k+5
\right\},
\\
LH
&=&
\left\{
i,j\in \mathbf{Z}
\colon
\
i< k-5,
\
k-3\leq j\leq k+3
\right\},
\\
HL
&=&
\left\{
i,j\in \mathbf{Z}
\colon
\
k-3\leq i\leq k+3,
\
j<k-5
\right\},
\\
HH
&=&
\left\{
i,j\in \mathbf{Z}
\colon
\
i,j>k+5,
\
|i-j|\leq 3
\right\}.
\end{eqnarray*}
We are going to estimate the  four terms 
in
(\ref{decomposition})
separately. For brevity  set
\begin{equation*}
\delta=\| V \|_{n/2\alpha}.
\end{equation*}
As mentioned above, we can
make
$\delta$
as small  as we wish by choosing
a small enough
$\rho$
in
(\ref{localised_equation}).
We will always assume that
$k$ is big enough, say
$k\geq 10$.

{\bf 3.}
By properties
of
$P_k$
and the Bernstein inequality
\begin{eqnarray*}
\|I\|_2
&
\lesssim
&
\sum_{i,j\in LL}
\|
P_{i}V
\,
P_{j}u
\|_2
\\
&
\lesssim
&
\sum_{i,j\in LL}
\| P_{i}V \|_\infty
\|P_{j} u \|_2
\\
&
\lesssim
&
2^{nk(2\alpha/n)}
\delta
\sum_{j=k-5}^{k+7}
\|P_{j} u \|_2.
\end{eqnarray*}
Term
$II$
is estimated exactly the same way.
It is convinient to record the final estimate in 
the following form
\begin{equation}
\label{I+II}
\|I\|_2
+
\|II\|_2
\lesssim
\delta
2^{  (2\alpha -s) k}
\sum_{j=k-5}^{k+7}
2^{s j}
\|P_j u\|_2
\end{equation}

{\bf 4.}
To estimate
$III$
we distinguish two cases.
First, assume that
\begin{equation}
\label{n_small}
n\leq 4\alpha,
\end{equation}
and hence
$n/2\alpha\leq 2$.
Apply the 
Holder inequality to derive 
\begin{eqnarray*}
\|III\|_2
&
\lesssim
&
\| P_{k-3\leq \cdot \leq k+3}V 
\,
P_{\cdot\leq 0} u\|_2
+
\sum_{j=1}^{k-5}
\| P_{k-3\leq \cdot\leq k+3}V 
\,
P_{j} u\|_2
\\
&
\lesssim
&
\| P_{k-3 \leq \cdot\leq k+3}V\|_2
\| P_{\cdot\leq 0} u\|_\infty
\\
&
&
+
\sum_{j=1}^{k-5}
\| P_{k-3\leq \cdot \leq k+3}V\|_2
\| P_{j} u\|_\infty
\\
&=&
X+Y.
\end{eqnarray*} 
From the Bernstein inequalities we  deduce that
\begin{eqnarray*}
X
&
\lesssim
&
2^
{
nk
(
(2\alpha/n)
-(1/2)
)
}
\|V\|_{n/2\alpha}
\| P_{\cdot\leq 0} u\|_2
\\
&
\lesssim
&
2^
{
2 \alpha k
-
(n/2)k
}
\delta
\| P_{\cdot\leq 0} u\|_2,
\end{eqnarray*}
and similarly
\begin{eqnarray*}
Y
&
\lesssim
&
\sum_{j=1}^{k-5}
2^
{
nk
(
(2\alpha/n) -(1/2)
)
}
\delta
\,
2^{nj/2}
\| P_{j} u\|_2
\\
&
\lesssim
&
\delta 
\sum_{j=1}^{k-5}
2^
{
2\alpha k -(n/2)k 
}
\,
2^{(n/2)j -sj}
\,
2^{sj}
\| P_{j} u\|_2.
\end{eqnarray*}
Consequently,
in the case of
(\ref{n_small}),
we can write the final estimate for
$III$
as
\begin{eqnarray}
\label{III_n_small}
\|III\|_2
&
\lesssim
&
\delta
2^
{
(2\alpha -(n/2))k
}
\|P_{\cdot\leq 0} u\|_2
\nonumber
\\
&
&
+
\delta 
2^{(2 \alpha -s) k}
\sum_{j=1}^{k-5}
\left(
2^{s j}
\| P_{j} u\|_2
\right)
2^
{
((n/2) - s )(j- k)
}.
\end{eqnarray}

Next assume that
\begin{equation}
\label{n_big}
n > 4\alpha.
\end{equation}
Hence
\begin{equation*}
\frac{2\alpha}{n}
+
\frac{n-4\alpha}{2n}
=
\frac{1}{2},
\quad
{\rm and}
\quad
\frac{n}{2\alpha},
\frac{2n}{n-4\alpha}
>
2
.
\end{equation*}
By the Holder inequality
\begin{eqnarray*}
\|III\|_2
&
\lesssim
&
\| P_{k-3\leq \cdot\leq k+3}V 
\|_{n/2\alpha}
\|
P_{\cdot\leq 0} u\|_{2n/(n-4\alpha)}
\\
&
&
+
\sum_{j=1}^{k-5}
\| 
P_{k-3\leq \cdot \leq k+3}V 
\|_{n/2\alpha}
\|
P_{j} u
\|_{2n/(n-4\alpha)}
\\
&
=
&
Z+W.
\end{eqnarray*}
The Bernstein inequalities
imply  that
\begin{equation*}
Z
\lesssim
\delta 
\|
P_{\cdot\leq 0} u
\|_2,
\end{equation*}
and
\begin{eqnarray*}
W
&
\lesssim
&
\delta
\sum_{j=1}^{k-5}
2^
{
nj
(
(1/2)
-(1/2)
+(2\alpha/n)
)
}
\|
P_{j} u
\|_2
\\
&
\lesssim
&
\delta
\sum_{j=1}^{k-5}
2^{(2 \alpha -s)j}
\,
2^{sj}
\|
P_{j} u
\|_2.
\end{eqnarray*}
Consequently, in the case of
(\ref{n_big}),
the final estimate for
$III$
can be written as
\begin{eqnarray}
\label{III_n_big}
\| III \|_2
&
\lesssim
&
\delta
\|P_{\cdot\leq 0} u\|_2
\nonumber
\\
&
&
+
\delta
2^{(2\alpha -s)k}
\sum_{j=1}^{k-5}
2^{(2\alpha -s )(j-k)}
\left(
2^{sj}
\|
P_{j} u
\|_2
\right)
.
\end{eqnarray}

{\bf 5.}
To estimate
$IV$
we 
also need to 
consider two cases.
First 
assume that
(\ref{n_small})
holds.
By  the Holder inequality
\begin{eqnarray*}
\|P_k(P_i V\, P_j u)\|_2
&
\lesssim
&
2^
{
nk/2
}
\|P_k(P_i V\, P_j u)\|_1
\\
&
\lesssim
&
2^
{
nk/2
}
\|P_i V\|_{n/2\alpha}
\| P_j u\|_{n/(n-2\alpha)}.
\end{eqnarray*}
According to 
(\ref{n_small})
we have
\begin{equation*}
\frac{n}{n-2\alpha}
\geq
2.
\end{equation*}
Therefore we can continue
with the help
of Bernstein inequality  and derive that
\begin{equation*}
\|P_k(P_i V\, P_j u)\|_2
\lesssim
2^{nk/2}
\delta
2^{nj((1/2) -1 +(2\alpha/n) )}
\|P_j u\|_2.
\end{equation*}
After the summation 
over 
$i$
and
$j$
lying in
the
$HH$
zone  we discover that
\begin{equation}
\label{IV_n_small}
\| IV \|_2
\lesssim
\delta 
2^{(2\alpha -s) k }
\sum_{j=k}^\infty
2^{((n/2) -2\alpha +s)(k-j)}
\left(
2^{s j}
\|P_j u\|_2
\right)
\end{equation}
provided
(\ref{n_small})
holds.

Next
assume that
(\ref{n_big})
holds.
Then
define
$q$,
$1\leq q \leq 2$
by writing
\begin{equation*}
\frac{1}{q}= \frac{1}{2}+ \frac{2\alpha}{n}.
\end{equation*}
Bernstein
and Holder
inequalities 
imply
that
\begin{eqnarray*}
\|P_k(P_i V\, P_j u)\|_2
&
\lesssim
&
2^
{
nk(2\alpha/n)
}
\|P_k(P_i V\, P_j u)\|_q
\\
&
\lesssim
&
2^
{
nk(2\alpha/n)
}
\|P_i V\, P_j u\|_q
\\
&
\lesssim
&
2^{2\alpha k}
\delta
\|P_j u\|_2.
\end{eqnarray*}
Summing this estimate over 
$i$
and 
$j$
in the
$HH$
region, 
we conclude 
that
in the case of
(\ref{n_big})
\begin{equation}
\label{IV_n_big}
\|IV\|_2
\lesssim
\delta
2^{(2\alpha-s)k}
\sum_{j=k}^\infty
2^{s(k-j)}
\left(
2^{s j}\|P_j u\|_2
\right)
.
\end{equation}

{\bf 6.}
Now we can prove the desired estimate
(\ref{goal_estimate}).
If 
(\ref{n_small})
holds, then substitute
(\ref{I+II}),
(\ref{III_n_small}),
and
(\ref{IV_n_small})
into
(\ref{decomposition}).
If
(\ref{n_big})
holds  then
use
(\ref{I+II}),
(\ref{III_n_big}),
and
(\ref{IV_n_big}).
To express  the result define
\begin{equation*}
\theta =
\left\{
\begin{array}{lcl}
\min\{1, (n/2) -s, (n/2) +s -2\alpha \}
&
{\rm for}
&
n\leq 4\alpha
\\
\min\{1, s, 2\alpha -s \}
&
{\rm for}
&
4\alpha < n.
\end{array}
\right.
\end{equation*}
According to assumptions of the theorem,
$\theta>0$.
Then we derive 
from
(\ref{first_estimate})
that for
$k\geq 10$
\begin{eqnarray}
\label{theta_decay}
2^{sk}
\|P_ku \|_2
&
\leq
&
{C_1(u, \rho)}{2^{-\theta k}}
\nonumber
\\
&
&
+
C_2(n,\alpha, s)
\delta
\sum_{j=0}^{\infty}
\left(
2^{s j}
\|P_j u\|_2
\right)
{2^{-\theta |j-k|}}.
\end{eqnarray}
For convenience set 
\begin{equation*}
a_k
=
2^{s k}
\|P_k u\|_2,
\quad
k=0,1,\ldots
.
\end{equation*}
We intend to use elemetary iteration 
Lemma~\ref{iter_lemma}
below
to bound the sequence
$\{a_k\}$.
First take
$\varepsilon =\theta/2$.
Next, find
$\rho>0$
such that in
(\ref{theta_decay})
we have
\begin{equation*}
\widetilde{\delta}
\stackrel{def}{=}
C_2 \delta
<
(
1-2^{-\varepsilon/100}
)
/2.
\end{equation*}
Then utilising
(\ref{theta_decay})
we can choose
$J=J(u,\rho)$
such that
\begin{equation*}
a_k
\leq
\frac{1}{2^{\varepsilon k}}
+
\widetilde{\delta}
\sum_{j=0}^\infty
\frac
{a_j}
{2^{2\varepsilon |k-j|}}
\quad
{\rm for}
\quad
k\geq J
\end{equation*}
with
$\widetilde{\delta}$ 
satisfying
(\ref{epsilon_small}).
Now, utilising the definition of
$H^s$-norm
find
$K=K(u, \rho)$
such that
\begin{equation*}
a_k
\leq
1
\quad
{\rm for}
\quad
k\geq K.
\end{equation*}
Finally set
\begin{equation*}
S=J+K.
\end{equation*}
All assumptions of  
Lemma~\ref{iter_lemma}
now hold and 
we derive
(\ref{goal_estimate}).
\end{proof}

The following lemma is a statement 
about number sequences. 
The proof of 
the lemma
is a 
careful but straightforward
iteration
of its  assumptions.
Actually we establish  a   stronger statement:
the proof shows that 
(\ref{sequence_decay})
holds even if in
(\ref{convolution})
we replace
$2\varepsilon$
by
any
$\varepsilon'>\varepsilon$.

\begin{lemma}
\label{iter_lemma}
Let
$\varepsilon>0$,
let
$\delta$
satisfy
\begin{equation}
\label{epsilon_small}
0< \delta < (1-2^{-\varepsilon})/2,
\end{equation} 
and let the sequence
$\{ a_k \}$
satisfy
\begin{eqnarray}
\label{convolution}
&&
0
\leq 
a_k 
\leq 
1
\quad
{\rm for}
\quad
k\geq S,
\nonumber
\\
&&
a_k
\leq
\frac{1}{2^{\varepsilon k}}
+
\delta
\sum_{j\geq 0}
\frac{a_j}{2^{2\varepsilon |k-j|}}
\quad
{\rm for}
\quad
k\geq S,
\end{eqnarray}
with some
$S\geq 0$.
Then 
\begin{equation}
\label{sequence_decay}
a_k
\leq
\frac{M}{2^{\varepsilon  k}},
\quad
k=0,1,\ldots,
\end{equation}
with a constant
$M\geq0$,
$M=M(\varepsilon, \delta, S, \| \{ a_k\}\|_{l^\infty}) $.
\end{lemma}

%\noindent
%{\it Proof.}
\begin{proof}
{\bf 1.}
From
the bounds on 
$a_k$
we derive at once that
\begin{equation}
\label{inductive_ass}
a_k
\leq
\frac{A}{2^{\varepsilon k}}
+
\delta
\sum_{j\geq S}
\frac{a_j}{2^{2\varepsilon |k-j|}}
\quad
{\rm for \quad all}
\quad
k\geq S
\end{equation}
with a constant
$A>0$,
$A=A(\varepsilon, \delta, S, \|\{a_k\}\|_{l^\infty})$.
Define
\begin{equation*}
C_\varepsilon =
2/
\left(
1-2^{-\varepsilon}
\right).
\end{equation*}
Then,  replacing
$a_j$
in
(\ref{inductive_ass})
by
$1$,
we also have 
\begin{eqnarray}
\label{induction_base}
a_k
&
\leq
&
\frac{A}{2^{\varepsilon k}}
+
\delta
\frac{2}{1-2^{-2\varepsilon}}
\nonumber
\\
&
\leq
&
\frac{A}{2^{\varepsilon k}}
+
\delta C_\varepsilon
\quad
{\rm for \quad all}
\quad
k\geq S.
\end{eqnarray}

{\bf 2.}
We claim that for any
$k\geq S$ 
and any
$N\geq 0$
the estimate
\begin{equation}
\label{N_iter_estimate}
a_k
\leq
\frac{A}{2^{\varepsilon k}}
\left(
1+\delta C_\varepsilon
+\cdots+
(\delta C_\varepsilon)^N
\right)
+
(\delta C_\varepsilon)^{N+1}
\end{equation}
holds.
Indeed,
for 
$N=0$
and all
$k\geq S$
this is just
(\ref{induction_base}).
Assume now that
(\ref{N_iter_estimate})
holds
for some
$N$
and all
$k\geq S$.
Then substitute   
(\ref{N_iter_estimate})
into
(\ref{inductive_ass})
to  discover that 
for any 
$k\geq S$
\begin{eqnarray*}
a_k
&
\leq
&
\frac{A}{2^{\varepsilon k}}
+
\delta A
\left(
1+\delta C_\varepsilon
+\cdots+
(\delta  C_\varepsilon)^{N}
\right)
\sum_{j\geq S}
\frac{1}{2^{\varepsilon j} 2^{2\varepsilon |j-k|}}
\\
&&
+\delta
(\delta C_\varepsilon)^{N+1}
\sum_{j\geq S}
\frac{1}{2^{2\varepsilon |j-k|}}
\\
&
\leq
&
\frac{A}{2^{\varepsilon k}}
\left(
1+\delta C_\varepsilon
+\cdots+
(\delta  C_\varepsilon)^{N+1}
\right)
+
(\delta  C_\varepsilon)^{N+2}
,
\end{eqnarray*}
because for 
$k\geq S$
\begin{eqnarray*}
\sum_{j\geq S}
\frac{1}{2^{\varepsilon j} 2^{2\varepsilon |j-k|}}
&
=
&
\sum_{j= S}^{k}
\frac{2^{\varepsilon  j}}{2^{2\varepsilon k}}
+
\sum_{j=k+1}^\infty
\frac{2^{2\varepsilon k}}{2^{3\varepsilon j}}
\\
&
\leq
&
\frac{1}{2^{\varepsilon k}}
\left(
\frac{1}{1-2^{-\varepsilon}}
\right)
+
\frac{1}{2^{\varepsilon k}}
\left(
\frac{1}{1-2^{-2\varepsilon}}
\right)
\\
&
\leq
&
\frac{C_\varepsilon }{2^{\varepsilon k}}.
\end{eqnarray*}
Hence
(\ref{N_iter_estimate})
is proved.

{\bf 3.}
Finally, sending 
$N$
to infinity in
(\ref{N_iter_estimate}),
we deduce according to
(\ref{epsilon_small})
that
$\delta C_\varepsilon<1$
and
\begin{equation*}
a_k\leq 
\left(
\frac{A}{1-\delta C_\varepsilon}
\right)
\frac{1}{2^{\varepsilon k}}
\quad
{\rm for \quad all}
\quad
k\geq S.
\end{equation*}
Thus
(\ref{sequence_decay})
holds.
\end{proof}
%$\Box$

%%%%%%%%%%%%%%%%%%%%%%%%%%%%%%%%%%%%%%%
%
%              BIBLIOGRAPHY
%
%%%%%%%%%%%%%%%%%%%%%%%%%%%%%%%%%%%%%%%

%
% BibTeX users please use
% \bibliographystyle{}
% \bibliography{}
%
% Non-BibTeX users please use

\end{document}